\documentclass{amsart}
\usepackage{amsmath, amssymb, euscript,  enumerate}

\newtheorem{thm}{Theorem}[section]
\newtheorem{cor}[thm]{Corollary}
\newtheorem{lem}[thm]{Lemma}
\newtheorem{prop}[thm]{Proposition}
\newtheorem{claim}[thm]{Claim}

\theoremstyle{remark}

\numberwithin{equation}{section}

\theoremstyle{definition}
\newtheorem{case}{Case}

\newcommand{\C}{\mathcal C}

\newcommand{\e}{\epsilon}
\newcommand{\dd}{\theta}

\newcommand{\ka}{\kappa}

\newcommand{\intd}{\int_0^{2\pi}}

\def\XXint#1#2#3{{\setbox0=\hbox{$#1{#2#3}{\intr}$}
     \vcenter{\hbox{$#2#3$}}\kern-.5\wd0}}

\begin{document}

\title{Classification of compact ancient solutions to the  curve
shortening flow}

\author{Panagiota Daskalopoulos$^*$}

\address{Department of
Mathematics, Columbia University, New York,
 USA}
\email{pdaskalo@math.columbia.edu}

\author{Richard Hamilton}

\address{Department of
Mathematics, Columbia University, New York,
 USA}
\email{hamilton@math.columbia.edu}

\author{Natasa Sesum$^{**}$}
\address{Department of Mathematics, Columbia University, New York,
USA}
\email{natasas@math.columbia.edu}

\thanks{$*:$ Partially supported
by NSF grant 0604657}

\begin{abstract}
We consider an embedded convex ancient solution $\Gamma_t$ to the
curve shortening flow in $\mathbb{R}^2$. We prove that there are only two
possibilities:  the family $\Gamma_t$ is either the family of
contracting circles,  which is a type I ancient solution,  or the family
of evolving Angenent ovals,  which correspond to a type II ancient
solution to the curve shortening flow.  We also give a necessary and sufficient 
curvature condition for an embedded, closed ancient solution to the curve shortening flow
to be convex.

\end{abstract}

\maketitle

\section{Introduction}

Let $M$ and $M'$ be two  Riemannian manifolds and $F: M\to M'$  a smooth
immersion. We can deform the immersion $F$ by the heat equation
$$\frac{\partial F}{\partial t} = \Delta F.$$
Since $\Delta F = \kappa \, \nu$, where $\kappa$ is the  trace of the
second fundamental form   and $\nu$ is the unit normal vector, we can
rewrite the previous evolution equation as
\begin{equation}
\label{eq-general}
\frac{\partial F}{\partial t} = \kappa\, \nu.
\end{equation}
The normal variation of the volume $V(M)$ is given by
$$\mathcal{D}V(M)(f) = -\int_M f\, \kappa \, dV$$
which implies that (\ref{eq-general}) is the gradient flow of the
volume functional.

We will  consider in this work  the special case when $M$ is an embedded curve
$\Gamma$ in $M'=\mathbb{R}^2$. Them, equation \eqref{eq-general} becomes the well studied
curve shortening flow. 
 In \cite{GH} Gage and Hamilton proved
that if $\Gamma$ is a convex curve embedded in $\mathbb{R}^2$, the
heat equation (\ref{eq-general}) shrinks $\Gamma_t$ to a point. In addition, the
curve remains convex and becomes asymptotically circular close to the
extinction time. In \cite{Gr} Grayson studied the evolution of non-convex embedded 
curves under (\ref{eq-general}). He proved that if $\Gamma$ is any embedded curve in $\mathbb{R}^2$, the solution
$\Gamma_t$ does not develop any singularities before it becomes strictly convex.

Let $\Gamma_t$ be an
embedded ancient solution to the curve shortening flow
(\ref{eq-general}).  If $s$ is the arclength along the curve, we can
express  the heat equation for the curve as a system
$$\frac{\partial x}{\partial t} = \frac{\partial^2 x}{\partial s^2}, \qquad  
\frac{\partial y}{\partial t} = \frac{\partial^2 y}{\partial s^2}.$$
The evolution for the curvature $\kappa$  of $\Gamma_t$ is given by
\begin{equation}
\label{eq-strictly-par}
\kappa_t  = \kappa_{ss} + \kappa^3
\end{equation}
which is a strictly parabolic equation. Let $\theta$ be the angle
between the tangent vector and the $x$ axis. For convex curves we can
use the angle $\theta$ as a parameter. It has been computed in  \cite{GH} that 
\begin{equation}\label{eqn-kappa}
\ka_t = \ka^2  \, \ka_{\dd\dd} + \ka^3.
\end{equation}
Equation 
It is also natural to look at the pressure function  $p=\ka^2$ which  evolves by
\begin{equation}\label{eqn-p}
p_t = p   \, p_{\dd\dd} - \frac 12 \,  p_\dd^2 +  2\, p^2.
\end{equation} 

It turns out that the evolution of the  family $\Gamma_t$ is completely described by the
evolution (\ref{eqn-kappa}) of the curvature $\kappa$. Gage and
Hamilton observed that a positive $2\pi$ periodic function
represents the curvature function of a simple closed strictly convex
$C^2$ plane curve if and only if
\begin{equation}
\label{eq-nec-cond}
\int_0^{2\pi}\frac{\cos\theta}{\kappa(\theta)}\, d\theta =
\int_0^{2\pi}\frac{\sin\theta}{\kappa(\theta)}\, d\theta = 0.
\end{equation}

\medskip

We will assume from now on that $\Gamma_t$ is an ancient solution of the curve shortening flow  defined on $(-\infty,T)$.  
We will also assume  that our extinction time $T = 0$. Let   $t_0 < 0$.
We define what it means to be a type I or type II ancient solution to 
(\ref{eqn-kappa}) as follows: 
\begin{itemize}
\item
It is type I if it satisfies $\sup_{\Gamma_t\times
(-\infty,t_0]}|t||p(x,t)| < \infty$.
\item
It is type II if $\sup_{\Gamma_t\times (-\infty,t_0]}|t||p(x,t)| = \infty$.
\end{itemize}

The  ancient solution to \eqref{eqn-kappa} defined  by $k(\theta,t) =
\frac{1}{\sqrt{-2t}}$ corresponds to a family of {\em contracting circles}. It is
clear that this solution is of type I and at the same time falls in a
category of contracting self-similar solutions (these are
solutions of the flow whose shapes change homothetically during the
evolution). In section \ref{sec-limit} we will prove the existence
of compact ancient solutions to (\ref{eqn-kappa}) that are not
self-similar. Since they have been discovered by Angenent we will
refer to them as to the {\em Angenent ovals}.

One very nice and important property of ancient solutions to the curve
shortening flow is that $\kappa_t \ge 0$. This fact follows from Hamilton's
Harnack estimate for convex curves (\cite{Ha}). By the strong maximum
principle, $\kappa(\cdot,t) > 0$ for all $t < 0$. If we start at any
time $t_0 \le 0$, Hamilton proved that
\begin{equation}
\label{eq-har-precise}
\kappa_t + \frac{\kappa}{2(t-t_0)} - \frac{\kappa_s^2}{k} \ge 0.
\end{equation}
Letting $t_0\to -\infty$ we get
\begin{equation}
\label{eq-harnack}
\kappa_t \ge 0.
\end{equation}

The main result of the paper is.

\begin{thm}
\label{thm-main}
Let $p(\theta,t) = \kappa^2(\theta,t)$ be an ancient solution
to (\ref{eqn-p}), definning a family of embedded closed convex curves in
$\mathbb{R}^2$ that evolve  by the curve shortening flow. Then 
\begin{itemize}
\item
either $p(\theta,t) = \frac{1}{(-2t)}$ which corresponds to a family of
contracting circles,
\item
or $p(\theta,t) = \lambda(\frac{1}{1- e^{2\lambda t}} - \sin^2(\theta + \gamma))$
which corresponds to Angenent ovals.
\end{itemize}
\end{thm}

One may ask whether the condition that the ancient solution be  convex is necessary 
in the previous classification  result. 
The following result gives a sufficient and necessary condition for a  embedded closed  ancient
solution to be convex. It is a consequence of the previous theorem.

\begin{cor}
\label{cor-equiv}
Let $\gamma(\cdot,t)$ be a  embedded closed  ancient solution of  the curve shortening flow (\ref{eq-general}). It is convex if and only if there are uniform constants $C_1, C_2$ so that $|\kappa(\cdot,t)| \le C_1$ and $\int_{\gamma_t}|\kappa|(s,t)\, ds \le C_2$ for all $t\in (-\infty,t_0)$, where $t_0 < 0$. 
\end{cor}

The organization of the paper is as follows:  In section
\ref{sec-lyapunov} we will introduce a monotone Lyapunov functional along the
flow that will help us classify the backward limits of our compact 
ancient solutions in section \ref{sec-limit}. It turns out that
the only non-trivial backward limits  are  grim reapers. In section
\ref{sec-classify} we will finish the proof of Theorem \ref{thm-main}
by analyzing the spectrum of the linearized  equation satisfied by  
the curvature function appropriately re-normalized.  In section \ref{sec-equiv} we give the proof of  Corollary \ref{cor-equiv}.

\section{The Lyapunov functional}
\label{sec-lyapunov}

Assume that we have an embedded  closed convex ancient solution of
the curve shortening flow (CSF).  By the result of M. Gage in \cite{G} we know that it collapses to a round point
at some time $T < \infty$. We will assume throughout this section that  $T = 0$.  
Let us consider the pressure function $p(\theta,t)$ defined by $p=\kappa^2$ which satisfies equation
\eqref{eqn-p} 

We define the Lyapunov  functional 
$$J(p) = \intd \left ( \frac{p_\dd^2}{p} - 4\, p \right )\, d\dd.$$
By direct computation we have:
\begin{equation}
\begin{split}
\frac {d}{d t} J(p(t)) &= \intd \left ( \frac{ 2 p_\dd p_{\dd t} }{p} - 
 \frac{ p_\dd^2  p_{t} }{p^2} - 4\, p_t \right )\, d \dd \nonumber \\
&=\intd \left ( \frac{ -2\,  p_{\dd\dd} p_{ t} }{p} + 
 \frac{ 2\, p_\dd^2  p_{t} }{p^2} -  \frac{  p_\dd^2  p_{t} }{p^2}-  4\, p_t \right )\, d \dd \nonumber \\
& = \intd \left ( \frac{ -2\, ( p   \, p_{\dd\dd} - \frac 12 \,  p_\dd^2 )\, p_t}{p^2} -  4\, p_t \right )\, d \dd \nonumber \\
& = \intd \left ( \frac{ -2\, ( p_t - 2\, p^2) \, p_t}{p^2} -  4\, p_t \right )\, d \dd \nonumber \\
&= - 2 \intd \frac{p_t^2}{p^2} \, d\dd.
\end{split} 
\end{equation}

\medskip

It is clear that $J(p(t)) \leq 0$ on the contracting circles and
the Angenent ovals. We will show next that this is true on any
ancient solution $p$.

\begin{lem} 
\label{lem-negative}
On an ancient solution $p$ of \eqref{eqn-p}
we  have $J(p(t)) \leq 0$, for all $t <0$. 
\end{lem}
\begin{proof} 

On an ancient solution we have $k_t \geq 0$ which gives $p_t \geq 0$. Hence, by
\eqref{eqn-p} we get
$$ p   \, p_{\dd\dd} - \frac 12 \,  p_\dd^2 +  2\, p^2 \geq 0.$$
Dividing by $p$ and integrating by parts gives:
$$\intd - \frac 12 \frac{p_\dd^2}{p} + 2\, p \, d\dd \geq 0$$
from which the inequality $J(p(t)) \leq 0$ readily follows. 
\end{proof}

\section{The limit of an ancient solution as $t\to -\infty$}
\label{sec-limit}

Our solution $p(\cdot,t) > 0$ being ancient has the property that $p_t
\ge 0$ (that follows  from (\ref{eq-harnack})), which implies the existence of the limiting function  
$$\tilde{p}(\theta) = \lim_{t\to  -\infty} p(\theta,t).$$  Using the
monotonicity of the Lyapunov functional $J(p)$, defined  in section
\ref{sec-lyapunov},  we will classify the  limiting function $\tilde p(\theta)$.  Before we prove the classification of
the limit we will show the existence of the Angenent ovals (ancient
compact solutions to the curve shortening flow that are not
self-similar).

\begin{lem}
\label{lem-ovals}
There exist  compact  ancient solutions to (\ref{eqn-p}) of the form
$$p(\theta,t) = \lambda(\frac{1}{1-e^{2\lambda t}} - \sin^2(\theta + \gamma))$$
where $\lambda > 0$ and $\gamma$ is a fixed angle.
\end{lem}

\begin{proof}
We  look for a solution to (\ref{eqn-p}) in the form $$p(\theta,t) =
a(t) - b(t)\sin^2(\theta + \gamma).$$  Then 
$$p_t = a'(t) - b'(t)\sin^2(\theta+\gamma), \quad 
p_{\theta} = -\sin 2(\theta+\gamma), \quad p_{\theta\theta} = -2\cos(\theta+\gamma).$$
If we plug everything in (\ref{eqn-p}) we get
$$a'(t) - b'(t)\sin^2(\theta+\gamma) = -2a(t)\, b(t) + 2a^2(t).$$
Since this has to hold for every $\theta$ we obtain that  $b'(t) = 0$,
which means $b(t) = C_0$. If we choose $\theta$ at which
$\sin(\theta + \gamma) = 0$, we have
$$a'(t) = -2C_0a(t) + 2a(t)^2.$$
Solving this ODE yields
$$a(t) = \frac{C_0}{1 - e^{2C_0(t - C_1)}}$$
for some constants $C_0, C_1$. We conclude that 
$$p(\theta,t) = C_0\, \left (\frac{1}{1 - e^{2C_0(t - C_1)}} - \sin^2(\theta + \gamma) \right ).$$ 
\end{proof} 

Observe that in the case of contracting circles the
$\lim_{t\to-\infty} p(\theta,t) = 0$ while  in the case of Angenent ovals
the $\lim_{t\to-\infty} p(\theta,t) = \lambda\cos^2(\theta + \gamma)$.
This motivates the following Proposition.

\begin{prop}
\label{prop-limit}
The limiting function $\tilde{p}(\theta)$ has always to be of the form 
$$\tilde p(\theta) = a\, \cos^2(\theta + b)$$ for some constants $a  \ge 0$ and $b$. Moreover, the convergence is smooth away from finitely many points 
in $[0,2\pi]$ which are exactly the zeros of $\tilde{p}(\theta)$.
\end{prop}

The proof of Proposition \ref{prop-limit} will be based on the Lyapunov functional identity
\begin{equation}
\label{eq-lyapunov}
\frac {d}{d t} J(p(t)) = - 2 \intd \frac{p_t^2}{p^2} \, d\dd
\end{equation}
shown in the previous section. We will prove  that the above identity implies  that the limit
$$\lim_{t \to -\infty} p_t =0$$
which shows  that $\tilde p$ must be  steady state of equation \eqref{eqn-p}. We will then derive that
 $\tilde p(\theta)=a\, \cos^2(\theta + b)$, for some   constants $a  \ge 0$ and $b$. To be able to pass to the limit in the equation \eqref{eqn-p}
 and the integral  identity \eqref{eq-lyapunov} we will first need to establish a priori derivative estimates  on the solution $p$.  

\begin{lem}
\label{lem-p-der}
There exists a uniform constant $C$ so that $p \le C$ and
$|p_{\theta}| \le C$ for all $t\le s_0 < 0$.
\end{lem}

\begin{proof}
By the Harnack estimate for ancient solutions to the curve shortening  flow we have
$p_t \ge 0$,  which implies $p(\theta,t) \le p(\theta,s_0) \le C$ for all $t\le s_0$.
Furthermore, $p_t\ge 0$ is equivalent to
\begin{equation}
\label{eq-equiv}
pp_{\theta\theta} - \frac{1}{2}p_{\theta}^2 + 2p^2 \ge 0.
\end{equation}
For each instant $t$ consider a point $\theta$ at which $p_{\theta}^2(\cdot,t)$ attains
its maximum. At that point $(p_{\theta}^2)_{\theta} = 0$, which implies $p_{\theta}p_{\theta\theta} = 0$.
If $p_{\theta} = 0$ at the maximum point of $p_{\theta}^2$, then $p_{\theta}\equiv 0$ for all
$\theta\in [0,2\pi]$ and the claim holds. If not, then $p_{\theta\theta} = 0$ and (\ref{eq-equiv})
implies that at the maximum point of $p_{\theta}^2$, we have 
$p_{\theta}^2 \le 4p^2 \le C$, for all $ t\le s_0$
and therefore the claim holds again.
\end{proof}

Recall that by  Lemma \ref{lem-negative} we have that 
\begin{equation}
\label{eq-weighted}
\int_0^{2\pi}\frac{p_{\theta}^2}{p}\, d\theta \le 4\int_0^{2\pi}p^2\, d\theta \le C,
\qquad  \mbox{for all} \,\, t\le s_0.
\end{equation}
Notice that (\ref{eq-equiv}) readily implies the bound from below 
\begin{equation}
\label{eq-lower-sec-der}
p_{\theta\theta} \ge -2p \ge -C, \quad  \mbox{for all} \,\,  t\le s_0.
\end{equation}

\begin{lem}
\label{lem-sec-der-L2}
There is a uniform constant $C$ so that
$$\int_0^{2\pi} p_{\theta\theta}^2\, d\theta \le C, \qquad \mbox{for all} 
\,\,  t  \le s_0-1.$$
\end{lem}

\begin{proof}
If we differentiate (\ref{eqn-p}) in $\theta$ we get
\begin{equation}
\label{eq-ev-first}
(p_{\theta})_t  = p(p_{\theta})_{\theta\theta} +
4pp_{\theta}
\end{equation}
which implies
$$\frac{1}{p}\, (p_{\theta}^2)_t  = p_{\theta}(p_{\theta})_{\theta\theta} + 4p_{\theta}^2$$
and therefore
$$\left ( \frac{p_{\theta}^2}{p} \right )_t = \frac{(p_{\theta}^2)_t}p  -
\frac{ p_\theta^2\,\,  p_t }{p^2}  \le  p_{\theta}(p_{\theta})_{\theta\theta} + 4p_{\theta}^2$$
where we use the fact that $p$ is an ancient solution implying that $p_t \ge 0$. Integrating  by parts in $\theta$ we get
$$\frac{d}{dt}\int_0^{2\pi}\frac{p_{\theta}^2}{p}\, d\theta \le -\int_0^{2\pi} p_{\theta\theta}^2\, d\theta +
4\int_0^{2\pi} p_{\theta}^2\, d\theta.$$
Using (\ref{eq-weighted}) and the fact that $|p_{\theta}| \le C$ uniformly in $t\le s_0$ we conclude that 
\begin{equation}
\label{eq-time-sec-der}
\int_{t-1}^{t+1}\int_0^{2\pi} p_{\theta\theta}^2\, d\theta ds \le C, \qquad  \mbox{for all} \,\, t\le s_0-1.
\end{equation}
Define $$I(t) = \int_0^{2\pi}p_{\theta\theta}^2(\theta,t)\, d\theta.$$
We will use the estimate (\ref{eq-time-sec-der}) to obtain a uniform 
bound on $I(t)$,   away from the extinction time.  If we differentiate (\ref{eq-ev-first}) in $\theta$ 
we get
$$(p_{\theta\theta})_t  = p\, (p_{\theta\theta})_{\theta\theta}
+ p_{\theta}\, p_{\theta\theta\theta} + 4\, p_{\theta}^2 +
4\,p\, p_{\theta\theta}.$$ 
If we multiply the previous equation by $p_{\theta\theta}$ and
integrate it over $\theta\in [0,2\pi]$ by parts, the uniform bounds on $p$ and $p_{\theta}$, the H\"older and the interpolation inequality yield
\begin{eqnarray}
\label{eq-third}
\frac{d}{dt}I(t) &=&
-2\int_0^{2\pi}p\,  p_{\theta\theta\theta}^2\, d\theta + 8\int_0^{2\pi} p_{\theta}^2\, p_{\theta\theta}\, d\theta
+ 8\int_0^{2\pi} p\,  p_{\theta\theta}^2\, d\theta \nonumber \\
&\le& C\, I(t) + C.
\end{eqnarray}
Let $\phi(s)$ be the cut off function such that $\phi(s) = 0$ for $s\in [t-1,t-1/2]$ and $\phi(s) = 1$ for
$s\in [t,t+1]$, $0 \le \phi(s) \le 1$ and $|\phi'(s)| \le 2$. Then,
$$\frac{d}{dt} (\phi(t)I(t)) \le C\phi(t) I(t) + C\phi(t) + \phi'(t) I(t).$$
Take any $s\in [t,t+1]$ and integrate the previous estimate over $[t-1,s]$ to get
$$I(s) \le C + C\int_{t-1}^{t+1} I(\tau)\, d\tau.$$
Since we can apply the previous analysis to every $t \le s_0-1$, 
by (\ref{eq-time-sec-der}) we get
$$I(t) \le C, \qquad  \mbox{for all} \,\, t\le s_0-1.$$
\end{proof}

As a  direct consequence of  the previous lemma and Sobolev's inequality we have:
 
\begin{cor}\label{cor-pdd}
We have $p_\dd(\cdot,t) \in C^{1/2}$ uniformly in time, i.e.
$$\|p_\dd(\cdot,t)\|_{C^{1/2}([0,2\pi])} \leq C, \qquad \mbox{for} \,\, t \leq s_0 -1.$$ 
\end{cor}

\medskip

We recall that since $p_t \ge 0$ and $p(\cdot,t) > 0$, the  pointwise limit $\tilde p=\lim_{t\to -\infty}p(\cdot,t) $ exists. 
Corollary \ref{cor-pdd} implies that 
$$p(\cdot,t) \stackrel{C^{1,1/2}}{\longrightarrow} \tilde{p}, \qquad  \mbox{as} \,\,\, t\to -\infty.$$ 

We will next classify the backwards limits $\tilde p$. We start with the following lemma which provides a bound from below on the distance between two zeros of   $\tilde p$.  This result follows as a consequence of
the:\\
{\bf Wirtinger's inequality (\cite{MD})}: if $f(a) = f(b) = 0$ with $b-a \le \lambda\pi$, then
$$\int_a^b f^2\, d\theta \le \lambda^2 \int_a^b f_{\theta}^2\, d\theta.$$

\begin{lem}
\label{lem-zeros}
There exists a $\delta > 0$ so that if $\tilde{p}(a) = \tilde{p}(b) = 0$, then either $|a-b| \ge \delta\pi$ or 
$\tilde{p}\equiv 0$ on $[a,b]$.
\end{lem}

\begin{proof}
Let $a < b\in [0,2\pi]$ be such that $\tilde{p}(a) = \tilde{p}(b) = 0$ and $b-a \le \delta\pi$.  
Wirtinger's inequality applied to our case yields
\begin{equation}
\label{eq-wirt}
\int_a^b\tilde{p}^2\, d\theta \le \e^2\int_a^b \tilde{p}_{\theta}^2\, d\theta.
\end{equation}
On our solution $p(\cdot,t)$ we have  $p_t \ge 0$, which implies
$$pp_{\theta\theta} - \frac{1}{2}p_{\theta}^2 + 2p^2 \ge 0.$$
If we integrate it over $[a,b]$ we get
$$\int_a^b p_{\theta}^2\, d\theta \le \frac{4}{3}\int_a^b p^2\, d\theta + 
\frac{3}{2}(p\, p_{\theta}(b,t) - p\, p_{\theta}(a,t)).$$
Letting $t\to -\infty$ in the previous inequality and  using that $p(\cdot,t) \stackrel{C^{1,1/2}}{\longrightarrow} \tilde{p}(\cdot)$ and that $\tilde{p}(a) = \tilde{p}(b) = 0$, we obtain 
\begin{equation}
\label{eq-har}
\int_a^b \tilde{p}_{\theta}^2\, d\theta \le \frac{4}{3}\int_a^b \tilde{p}^2\, d\theta.
\end{equation}
Combining (\ref{eq-wirt}) and (\ref{eq-har}) yields a contradiction if  $\delta^2 < \frac{3}{4}$,  unless $\tilde{p} \equiv 0$ on $[a,b]$.
\end{proof}

Take the same $\delta > 0$ as in Lemma \ref{lem-zeros}. A simple corollary of it is the following observation.

\begin{cor}
\label{cor-poss}
There are finitely many points $0 \le \theta_1 < \dots \theta_n \le 2\pi$ so that $\tilde{p}(\theta_i) =0$ and
\begin{itemize}
\item
either $\tilde{p}(\cdot) \equiv 0$ on $[\theta_{k-1},\theta_k]$ 
\item
or $\tilde{p} > 0$ on $[\theta_{k-1},\theta_k]$ and $\theta_k - \theta_{k-1} \ge \delta\pi$.
\end{itemize}
\end{cor}

We will next prove: 

\begin{lem}
\label{lem-form-interval}
The solution $\tilde{p}$ has to be of the form $\tilde p(\theta)=a\cos^2(\theta + b)$,  for some  constants $a,b$
on each of the intervals $[\theta_{k-1},\theta_k]$ discussed  above.
\end{lem}

\begin{proof}
If we take an interval $[\theta_{k-1},\theta_k]$ such as in the
first case of Corollary \ref{cor-poss}, we can just take $a = 0$ and
$b$ arbitrary. Lets analyze the form of our solution $\tilde{p}$ on
an interval $[\theta_{k-1},\theta_k]$ given by the second
possibility discussed in Corollary \ref{cor-poss}. On any open
subset $I$, compactly contained in $[\theta_{k-1},\theta_k]$, we
have that $\tilde{p} \ge \eta > 0$ for some constant $\eta$ (that
depends on $I$). Since $p_t \ge 0$, we have $p(\theta,t) \ge \eta$ for all
$(\theta,t) \in I\times (-\infty,0)$. This implies the uniform in
$t$ parabolicity of (\ref{eqn-p}) on $I$. Standard local parabolic
estimates give us uniform estimates on higher order derivatives of
$p(\cdot,t)$, besides the uniform $C^{1,1/2}$ and $W^{2,2}$
estimates provided already by Lemma \ref{lem-sec-der-L2} and
Corollary \ref{cor-pdd}. 

\begin{claim}
Our limit $\tilde{p}(\theta)$ is a steady state of equation (\ref{eqn-p}),
that is, it satisfies
\begin{equation}
\label{eq-steady-state}
\tilde{p}\tilde{p}_{\theta\theta} - \frac{1}{2}\tilde{p}_{\theta}^2 + 2\tilde{p}^2 = 0
\end{equation}
on $(\theta_{k-1},\theta_k)$ in the classical strong sense.
\end{claim}

\begin{proof}
By Lemma \ref{lem-negative}, $J(p(t)) \le 0$. On the other hand, we have
$$J(p(t)) \ge -4\int_0^{2\pi} p\, d\theta \ge -C$$ by Lemma \ref{lem-p-der}.
Since $J(p(t))$ is monotone, it follows there exists a finite limit
$$J(p(-\infty)) = \lim_{t\to -\infty} J(p(t)).$$
Take any $A > 0$. Integrating (\ref{eq-lyapunov}) in $t$ yields
\begin{eqnarray*}
J(p(t+A)) - J(p(t)) &=& -\int_t^{t+A}\int_0^{2\pi}\frac{p_t^2}{p^2}\, d\theta\, ds \\
&=& -\int_0^A\int_0^{2\pi}\frac{p_t^2}{p^2}(t+s)\, d\theta\, ds.
\end{eqnarray*}
Observe that $\lim_{t\to -\infty} J(p(t)) = \lim_{t\to -\infty}
J(p(t+A))$ and therefore, we have 
$$\lim_{t\to -\infty}\int_0^A\int_0^{2\pi}\frac{p_t^2}{p^2}(t+s)\, d\theta\, ds = 0.$$
The uniform higher order derivative estimates away from the zeros of $\tilde{p}$,
Arzela-Ascoli theorem, the uniqueness of the limit $\tilde{p}$ and
Fatou's lemma now imply that 
$$\int_0^A\int_{\theta_{k-1}}^{\theta_k}\frac{\tilde{p}_t^2}{\tilde{p}^2}\, d\theta\, ds 
\le \lim_{t\to -\infty}\int_0^A\int_{\theta_{k-1}}^{\theta_k}\frac{p_t^2}{p^2}(t+s)\, 
d\theta\, ds = 0$$
which yields to the equality $\tilde{p}_t = 0$.  Hence, 
$$\tilde{p}\tilde{p}_{\theta\theta} - \frac{1}{2}\tilde{p}_{\theta}^2 + 2\tilde{p}^2 = 0$$
on $(\theta_{k-1},\theta_k)$ in the classical  strong sense.
\end{proof}

By using the equation (\ref{eq-steady-state}) we compute
\begin{eqnarray*}
 \left (\frac{\tilde{p}_{\theta}^2}{\tilde{p}} \right )_{\theta} &=& 
\frac{2\tilde{p}_{\theta}\tilde{p}_{\theta\theta}\tilde{p} -\tilde{p}_{\theta}^3}{\tilde{p}^2} \\
 &=& \frac{2\tilde{p}_{\theta}(\tilde p \, \tilde{p}_{\theta\theta} - \frac{1}{2}\tilde{p}_{\theta}^2)}{\tilde{p}^2} =
 \frac{2\, \tilde{p}_{\theta}}{\tilde{p}^2}\, (-2\tilde{p}^2) \\
 &=& -4\, \tilde{p}_{\theta}
\end{eqnarray*} 
which, after integration, implies the differential  equation
$$\frac{\tilde{p}_{\theta}^2}{\tilde{p}} + 4\tilde{p} =
C.$$  Solving this ODE yields to explicit solutions of the form
$\tilde{p}(\theta) = a\cos^2(\theta + b)$ for some constants $a,b$, 
where $a \ge 0$ since we are looking for nonnegative solutions.
\end{proof}

The previous lemma shows that  there are finitely many  $0 \le
\theta_1 \le \dots \le \theta_n \le 2\pi$ so that 
\begin{equation}\label{solutions}
\tilde{p}(\theta_k)
= 0 \qquad \mbox{and} \qquad  \tilde{p}(\theta) = a_k\cos^2(\theta+b_k),  \quad \mbox{for} \,\,  \theta\in
[\theta_{k-1},\theta_k].
\end{equation}
In order to conclude the proof of 
our Proposition \ref{prop-limit}, we  will  show next that on the whole interval
$[0,2\pi]$, we have $\tilde{p}(\theta) = a\cos^2(\theta + b)$, for some unique
constants $a \ge 0$ and $b$.

\begin{lem}
\label{lem-form}
There are unique constants $a > 0$ and $b$ so that $\tilde{p}(\theta) = a\cos^2(\theta + b)$,  for all
$\theta\in [0,2\pi]$.
\end{lem}

\begin{proof}
Recall that \eqref{solutions} hold  and observe that at $\theta_k$ we have $a_k\cos^2(\theta_k + b_k) =
a_{k+1}\cos^2(\theta_k + b_{k+1}) = 0$. This implies the following  three  possibilities:
\begin{itemize}
\item
$a_k$ or $a_{k+1} = 0$; 
\item
both $\theta_k + b_k$ and $\theta_k + b_{k+1}$ are at the same time $\pi/2$ or $3\pi/2$,
which implies $b_k = b_{k+1}$;  or
\item
$|b_k - b_{k+1}| = \pi$.
\end{itemize}
In the last two cases,  we have $\cos^2(\theta + b_k) = \cos^2(\theta + b_{k+1})$ for
$\theta\in [0,2\pi]$. It follows  that our solution $\tilde{p}(\theta)$
can be written in the form 
\begin{equation}
\label{eq-form-p}
\tilde{p}(\theta) = a_k\cos^2(\theta + b) \qquad  \mbox{for} \,\, \theta\in [\theta_{k-1},\theta_k].
\end{equation}
The definition of $\tilde{p}$ can change only at its zeros and from
(\ref{eq-form-p}) we see that the only possibilities are when $\theta
= \frac{(2k+1)\pi}{2} - b$ where $k\in \mathbb{Z}$.  Because of the
periodicity and the fact that $\cos^2(\theta + b) = \cos^2(\theta + b
- \pi)$, we may take $b \in [-\pi,0]$. Since $\theta \in [0,2\pi]$,
depending on $b$, we have that either $\tilde{p}(\theta)$ is zero at
$0, \pi, 2\pi$, in which case $b = -\frac{\pi}{2}$, 
or there are two zeros $0 < \theta_1 < \theta_2 < 2\pi$
so that $\theta_2 = \theta_1 + \pi$. In the former case 
\begin{equation}  
\tilde{p}(\theta) =
\begin{cases}
a_1\cos^2\theta, \quad &\mbox{if} \,\,\, \theta\in [0,\pi] \\
a_2\cos^2\theta, \quad &\mbox{if} \,\,\, \theta\in [\pi,2\pi].
\end{cases}
\end{equation} 
Since $\tilde{p}(\theta)$ is a $2\pi$-periodic function, we have
$a_1 = \tilde{p}(0) = \tilde{p}(2\pi) = a_2$ and we are done.
Lets consider the case when $\tilde{p}$ has two zeros 
$0 < \theta_1 < \theta_2 < 2\pi$. Assume
\begin{equation}  
\tilde{p}(\theta) =
\begin{cases}
a_1\cos^2(\theta+b), \quad &\mbox{if} \,\,\, \theta\in [0,\theta_1] \\
a_2\cos^2(\theta+b), \quad &\mbox{if} \,\,\, \theta\in [\theta_1,\theta_2] \\
a_3\cos^2(\theta+b), \qquad &\mbox{if} \,\,\, \theta \in [\theta_2,2\pi].
\end{cases}
\end{equation}
Since $a_1\cos^2 b = \tilde{p}(0) = \tilde{p}(2\pi) = a_3\cos^2 b$ and $\cos b \neq 0$ we have
$a_1 = a_3$. By (\ref{eq-nec-cond}) we have
$$\int_0^{2\pi}\frac{\cos\theta}{\tilde{\kappa}}\, d\theta = \int_0^{2\pi}
\frac{\sin\theta}{\tilde{\kappa}}\, d\theta = 0$$
and therefore,
$$\int_0^{2\pi}\frac{\cos(\theta + b)}{\tilde{\kappa}}\, d\theta = 0.$$
This implies that 
\begin{eqnarray*}
0 &=& \int_0^{\theta_1}\frac{\cos(\theta + b)}{\tilde{\kappa}}\, d\theta -
\int_{\theta_1}^{\theta_1 + \pi}\frac{\cos(\theta + b)}{\tilde{\kappa}}\, d\theta
+ \int_{\theta_1 + \pi}^{2\pi} \frac{\cos(\theta + b)}{\tilde{\kappa}}\, d\theta \\
&=& \frac{\theta_1}{\sqrt{a_1}} - \frac{\pi}{\sqrt{a_2}} + \frac{\pi - \theta_1}{\sqrt{a_1}}
= \pi(\frac{1}{\sqrt{a_1}} - \frac{1}{\sqrt{a_2}})
\end{eqnarray*}
which proves that  $a_1 = a_2$, finishing  the proof of Lemma \ref{lem-form}
\end{proof} 

\section{The classification of ancient solutions to (\ref{eqn-p})}
\label{sec-classify}

In this section we will prove that the only closed convex ancient solutions to the
curve shortening flow, whose extinction time is given by $t = 0$, are
the contracting circles and Angenent ovals.  More precisely, we will
prove the following theorem.

\begin{thm}
\label{thm-classification}
The only ancient solutions to (\ref{eqn-p}), corresponding to closed convex curves 
in $\mathbb{R}^2$ evolving by the curve shortening flow are:
\begin{enumerate}
\item[(i)]
either $p(\theta,t) = \frac{1}{(-2t)}$, which corresponds to contracting circles,
or 
\item[(ii)]
$p(\theta,t) = \lambda(\frac{1}{1-e^{2\lambda t}} - sin^2\theta)$, for
a parameter $\lambda > 0$, which corresponds to the Angenent ovals.
\end{enumerate}
\end{thm}

We will prove the theorem  by introducing another monotone functional along the flow. To analyze the behavior of our functional as $t \to -\infty$, we will   use  the results of the previous section. We will also need to analyze   behavior of  the
functional as $t\to 0$, which   will use the following result of
Gage and Hamilton (\cite{GH}).

\begin{thm}[Gage, Hamilton]
\label{thm-gh}
If $\Gamma$ is a closed convex curve embedded in the plane
$\mathbb{R}^2$, the curve shortening flow shrinks $\Gamma$ to a point
in a circular manner. Moreover, the curvature and all its derivatives
of a rescaled curve shortening flow converge exponentially to $1$ and
$0$ respectively, with the rate $e^{-2\eta \tau}$, where $\tau$ is our
new time variable introduced below and $\eta$ is any constant in
$(0,1)$.
\end{thm}

Let us remark here that in order to  determine the rate of convergence of the curvature $\kappa$ and its
derivatives in \cite{GH}, it turned out to be the most productive to consider the
evolution of the normalized curvature $\tilde{\kappa}$,  where the
normalization is chosen so that the related convex curve encloses an
area $\pi$. The rescaled curvature   $\tilde{\kappa}$ is defined by
\begin{equation}\label{tilde-k}
\tilde{\kappa}(\theta,\tau) = k(\theta,t)\sqrt{-2t}, \qquad \mbox{with} \,\, \tau = -\frac{1}{2}\log(-t).
\end{equation} The evolution equation for
$\tilde{\kappa}$ is
\begin{equation}
\label{eq-tilde-k}
\tilde{\kappa}_t = \tilde{\kappa}^2\tilde{\kappa}_{\theta\theta} + \tilde{\kappa}^3 - \tilde{\kappa}.
\end{equation}

\medskip

We will now define our new monotone functional. Denote by
$$\alpha(\theta,t) := p_{\theta}(\theta,t).$$  By using (\ref{eqn-p}) it easily follows that
\begin{equation}
\label{eq-der-p}
\alpha_t  = p\, (\alpha_{\theta\theta} + 4\, \alpha).
\end{equation}
We introduce the functional
$$I(\alpha) = \int_0^{2\pi}(\alpha_{\theta}^2 - 4\alpha^2)\, d\theta.$$
The following lemma shows the monotonicity of $T(\alpha)$ in time. 

\begin{lem}
\label{lem-mon-funct}
$I(\alpha(t))$ is decreasing along the flow (\ref{eq-der-p}). Moreover,
$$\frac{d}{dt} I(\alpha(t)) =  -2\int_0^{2\pi}\frac{\alpha_t^2}{p}\, d\theta.$$
\end{lem}

\begin{proof}
We compute
\begin{eqnarray*}
\frac{d}{dt} I(\alpha(t)) &=& \int_0^{2\pi} (2\alpha_{\theta}\alpha_{\theta t} - 8\alpha\alpha_t)\, d\theta 
\\
&=& -2\int_0^{2\pi}2\alpha_{\theta\theta}\alpha_t\, d\theta - 8\int_0^{2\pi}\alpha\alpha_t\, d\theta \\
&=& \int_0^{2\pi}\frac{2(\alpha_t - 4\alpha p)\alpha_t}{p}\, d\theta - 8\int_0^{2\pi}\alpha\alpha_t\, d\theta \\
&=& -2\int_0^{2\pi}\frac{\alpha_t^2}{p}\, d\theta.
\end{eqnarray*}
\end{proof}

An easy computation shows that $I(\alpha(t)) \equiv 0$ on both  the
circles,  which are defined by $p(\theta,t) = \frac{1}{(-2t)}$, and on the 
Angenent ovals which are given by $p(\theta,t) = \lambda(\frac{1}{1-e^{2\lambda
t}} - \sin^2 \theta)$, for $\lambda > 0$. This motivates  the following proposition.

\begin{prop}
\label{prop-ancient}
For any ancient solution to (\ref{eqn-p}), $I(\alpha(t)) \equiv 0$ for all $t\in (-\infty,0)$.
\end{prop}

\begin{proof}
By Lemma \ref{lem-mon-funct}, $I(\alpha(t))$ is decreasing in time.
By Proposition \ref{prop-limit} we know that $p(\theta,t) \to
a\, \cos^2(\theta + b) =:\tilde{p}$ in $C^{1,1/2}$ norm. Assume with no
loss of generality that $a = 1$ and $b = 0$. SInce $p_t \ge 0$, we have
$p(\theta,t) \ge \cos^2\theta$ for all $t < 0$.

\begin{lem}
\label{lem-conv-sec}
There is a sequence $t_i\to -\infty$ so that 
$$\lim_{i\to \infty}\int_0^{2\pi}p_{\theta\theta}^2\, d\theta = \int_0^{2\pi}\tilde{p}_{\theta\theta}^2\, d\theta.$$
Furthermore, the $\lim_{t\to -\infty} I(\alpha(t)) = 0$.
\end{lem}

\begin{proof}
By (\ref{eq-har-precise}) we have
$$pp_{\theta\theta} - \frac{1}{2}p_{\theta}^2 + 2p^2 \ge \frac{p_{\theta}^2}{2},$$
which implies
\begin{equation}
\label{eq-har-p}
\frac{p_{\theta}^2}{p} \le p_{\theta\theta} + 2p.
\end{equation}
At the maximum of $\frac{p_{\theta}^2}{p}$ we have $p_{\theta\theta} = \frac{p_{\theta}^2}{2p}$, 
which together with (\ref{eq-har-p}) implies
\begin{equation}
\label{eq-p-quot}
\frac{p_{\theta}^2}{p} \le C.
\end{equation}
By Lemma \ref{lem-sec-der-L2} and (\ref{eq-third}) we have
$$\int_{t-1}^{t+1}\int_0^{2\pi} pp_{\theta\theta\theta}^2\, d\theta\, ds \le C, \qquad \mbox{for all}
\,\,\, t \le s_0-1.$$
This implies there is a subsequence $t_i\to -\infty$ so that 
\begin{equation}
\label{eq-seq-third}
\int_0^{2\pi}pp_{\theta\theta\theta}^2(t_i)\, d\theta \le C,
\end{equation}
for a uniform constant $C$. Let $0 < \beta < \frac{1}{2}$. Then for $t = t_i$,
using integration by parts, (\ref{eq-seq-third}) and H\"older inequality,
\begin{eqnarray}
\label{eq-p10}
|\int_0^{2\pi} (p_{\theta\theta}^2 - \tilde{p}_{\theta\theta}^2)|\, d\theta &=&
|\int_0^{2\pi}(p_{\theta\theta} - \tilde{p}_{\theta\theta})\cdot (p_{\theta\theta} + \tilde{p}_{\theta\theta})\, d\theta| \nonumber \\
&=& |\int_0^{2\pi} (p_{\theta} - \tilde{p}_{\theta})\cdot(p_{\theta\theta\theta} - \tilde{p}_{\theta\theta\theta})\, d\theta| \nonumber\\
&\le& 2(\int_0^{2\pi}\frac{|p_{\theta}-\tilde{p}_{\theta}|^2}{p}\, d\theta)^{1/2}\cdot(\int_0^{2\pi}p\cdot(p_{\theta\theta\theta}^2
+ \tilde{p}_{\theta\theta\theta}^2)\, d\theta)^{1/2} \nonumber \\
&\le& C(\int_0^{2\pi}\frac{|p_{\theta}-\tilde{p}_{\theta}|^2}{p}\, d\theta)^{1/2}
\end{eqnarray}
Furthermore,
\begin{eqnarray}
\label{eq-p20}
\int_0^{2\pi}\frac{|p_{\theta} - \tilde{p}_{\theta}|^2}{p}\, d\theta &=&
\int_0^{2\pi}|p_{\theta} - \tilde{p}_{\theta}|^{\beta}\cdot \frac{|p_{\theta}-\tilde{p}_{\theta}|^{2-\beta}}
{p^{\frac{2-\beta}{3}}}\cdot\frac{1}{p^{\frac{1+\beta}{3}}} \nonumber \\
&\le& C\sup_{\theta\in [0,2\pi]}|p_{\theta} - \tilde{p}_{\theta}|^{\beta}\cdot\int_0^{2\pi}\frac{d\theta}{p^{\frac{1+\beta}{3}}} \nonumber \\
&\le& C\sup_{\theta\in [0,2\pi]}|p_{\theta} - \tilde{p}_{\theta}|^{\beta}\cdot
\int_0^{2\pi}\frac{d\theta}{(\cos\theta)^{\frac{2+2\beta}{3}}},
\end{eqnarray}
where we have used that $p(\theta,t) \ge \cos^2\theta$ and that 
\begin{eqnarray*}
\frac{|p_{\theta}-\tilde{p}_{\theta}|}{p^{1/3}} &\le& p^{1/6}\cdot (\frac{|p_{\theta}|}{p^{1/2}} +
\frac{|\tilde{p}_{\theta}|}{p^{1/2}}) \\
&\le& C(\frac{|p_{\theta}|}{p^{1/2}} + \frac{|\tilde{p}_{\theta}|}{\tilde{p}^{1/2}}) \\
&\le& \tilde{C},
\end{eqnarray*} 
by (\ref{eq-p-quot}). Since 
$$\int_0^{2\pi}\frac{d\theta}{(\cos\theta)^{\frac{2+2\beta}{3}}} < \infty \qquad \mbox{for} \,\,\, 
\beta < 1/2,$$
by (\ref{eq-p10}), (\ref{eq-p20}) and by the fact that $p(\cdot,t) \stackrel{\C^{1,1/2}}{\longrightarrow} \tilde{p}(\cdot,t)$
as $t\to -\infty$ we have
$$\lim_{i\to\infty}\int_0^{2\pi}\alpha_{\theta}^2\, d\theta = \int_0^{2\pi}\tilde{p}_{\theta\theta}^2\, d\theta
= 4\int_0^{2\pi}\cos^2\theta\, d\theta.$$
We also have $\alpha(\theta,t) \to -\, \sin 2\theta$ uniformly in
$\theta\in [0,2\pi]$ as $t\to -\infty$. Combining that and Lemma \ref{lem-conv-sec} yields that
\begin{equation}
\label{eq-lim-inf}
\lim_{t_i\to -\infty} \int_0^{2\pi} (\alpha_{\theta}^2 - 4\alpha^2)(t_i)\, d\theta =
4 \int_{0}^{2\pi}(\cos^2 \theta - \sin^2 \theta)\, d \theta = 0.
\end{equation}
Since $I(\alpha(t))$ is monotone decreasing along the flow and since (\ref{eq-lim-inf}), we also have
$$\lim_{t\to -\infty}I(\alpha(t)) = 0.$$
\end{proof}

To conclude the proposition it is enough to show that the $\lim_{t\to 0} I(\alpha(t)) = 0$ since that combined with the monotonicity of 
$I(\alpha(t))$ and (\ref{eq-lim-inf}) would yield that $I(\alpha(t)) \equiv 0$ for all $t\in (-\infty,0)$. This is shown in the next lemma. \end{proof}

\begin{lem}
\label{lem-I-zero}
$\lim_{t\to 0} I(\alpha(t)) = 0$.
\end{lem}

\begin{proof}
To prove the lemma we will analyze the normalized flow
(\ref{eq-tilde-k}) since due to \cite{GH} we have good decay
estimates on $\tilde{\kappa} - 1$ and the derivatives of $\tilde{\kappa}$.
Notice that when $\tilde \kappa$ is defined by eqref{tilde-k}, then
$$I(\alpha(t)) = e^{4\tau}\int_0^{2\pi}(\tilde{\alpha}_{\theta}^2 - 4\tilde{\alpha}^2)\, d\theta$$
with $\tilde \alpha$ denoting the corresponding rescaled $\alpha = p_\theta$. 
By \cite{GH} we have
\begin{equation}
\label{eq-decay}
|\tilde{\kappa} - 1| \le C(\eta) e^{-2\eta\tau}, \qquad  |\frac{\partial^m}{\partial \theta^m}\tilde{\kappa}| \le C_m(\eta)e^{-2\eta\tau}
\end{equation}
where $\eta\in (0,1)$. These estimates imply the bounds 
$$|\tilde{\alpha}| \le Ce^{-2\eta\tau}, \qquad |\tilde{\alpha}_{\theta}| \le Ce^{-2\eta\tau}.$$
We will utilize those bounds to conclude that  $\lim_{t\to 0} I(\alpha(t)) = 0$. 
In order to show that, we will analyze the linearization of (\ref{eq-tilde-k}) around $\tilde \kappa = 1$.
It is easy to see that
$$\frac{\partial}{\partial \tau}(\tilde{\kappa} - 1) = (\tilde{\kappa}-1)_{\theta\theta} + 2(\tilde{\kappa}-1) +
(\tilde{\kappa}-1)(\tilde{\kappa}-1)_{\theta\theta}(\tilde{\kappa}+1) + (\tilde{\kappa}-1)^2(\tilde{\kappa}+2)$$ 
which we can rewrite as
$$\frac{\partial}{\partial \tau}(\tilde{\kappa}-1) = \mathcal{L}(\tilde{\kappa}-1) + \mathcal{R}(\tilde{\kappa}-1)$$
where 
$$\mathcal{L}(f) = f_{\theta\theta} + 2f \qquad  \mbox{and} \qquad 
\mathcal{R} = ff_{\theta\theta}(f+2) + f^2(f+3).$$ 
Note that $\mathcal{R}(\tilde{\kappa}-1)$ is an error term that is quadratic in
$\tilde{\kappa}-1$ and its derivatives which all converge to zero
exponentially as $\tau\to \infty$.  The spectrum for $\mathcal{L}$ on
an interval $[0,2\pi]$ is given by
$$\lambda_l = 2 - l^2, \qquad l\ge 0$$
with the corresponding eigenvectors 
$$f_l(\theta) = \cos(l\theta) \qquad  \mbox{and} \qquad  g_l(\theta) = \sin(l\theta).$$
Denote by $w = \tilde{\kappa} - 1$. The semigroup representation formula for $w$ gives
\begin{equation}
\label{eq-semi-rep}
w(\theta,s) = e^{s\mathcal{L}}w(\theta,0) + 
\int_0^s e^{(s-\tau)\mathcal{L}}\, \mathcal{R}(w(\theta,\tau))\, d\tau. 
\end{equation}
In the space of continuous functions on $[0,2\pi]$ we can look at the
trigonometric system of functions $\{\cos(l\theta)\}_{l \ge
0}$ and $\{\sin(l \theta)\}_{l \ge 0}$. It is easy to check
that this is an orthogonal basis in $C_0[0,2\pi]$ with respect to the
inner product given by
$$(f,g) = \int_0^{2\pi} f(\theta)\, g(\theta)\, d\theta.$$
Hence  
$$w(\theta,0) = \sum_{l\ge 0} (\alpha_l f_l + \beta_l g_l).$$           
By definition, we have 
\begin{equation}
\begin{split}
\label{eq-zero-part}
e^{s\mathcal{L}}w(\theta,0) &= \alpha_0 e^{2s} + (\alpha_1\cos\theta + 
\beta_1\sin\theta) e^s \\&+ (\alpha_2\cos 2\theta + \beta_2\sin 2\theta) e^{-2s} + o(e^{-2s}). 
\end{split}
\end{equation}
Similarly, we have 
\begin{equation}
\label{eq-rep-R}
\mathcal{R}(w(\theta,\tau)) = \sum_{l\ge 0} (\alpha_l(\tau) f_l + \beta_l(\tau)g_l)
\end{equation}
which implies that 
\begin{equation}
\label{eq-error-part}
\begin{split}
\int_0^s e^{(s-\tau)\mathcal{L}} \mathcal{R}(w(\theta,\tau)) &= \int_0^s \left ( \alpha_0(\tau)\, e^{2(s-\tau)} 
+ [\, \alpha_1(\tau)\cos\theta + \beta_1(\tau)\sin\theta \, ]\, e^{(s-\tau)} \right . \\
&+ \left .  [\, \alpha_2(\tau)\cos 2\theta + \beta_2(\tau)\sin 2\theta \, ] \, e^{-2(s-\tau)} \,
\right ) \, d\tau 
\\& + \int_0^s \left (  \sum_{l\ge 3}   \, [\, \alpha_l(\tau)\cos l\theta + \beta_l(\tau)\sin l\theta\,  ] \, e^{\lambda_l(s-\tau)} \right ) \, d\tau
\end{split}
\end{equation}
where $\lambda_l < -2$ for $l\ge 3$.

\begin{claim} 
\label{claim-error} The lest term above can be estimated as: 
$$\left |\int_0^s \sum_{l\ge 3} [ \alpha_l(\tau)\cos\frac{l\theta}{2}  + \beta_l(\tau)\sin\frac{l\theta}{2})] \,e^{\lambda_l(s-\tau)}\,  d\tau \right | = o(e^{-2s}), \quad \mbox{as} \,\, s\to\infty$$
\end{claim}

\begin{proof}
We have $|\mathcal{R}(w(\theta,\tau))| \le Ce^{-4\eta \tau}$, which
implies $||\mathcal{R}(w(\cdot,\tau)||_{L^2[0,2\pi]} \le C
e^{-4\eta\tau}$. Since $\{\cos(l\theta),
\sin(l \theta)\}_{l\ge 0}$ is an orthogonal trigonometric
system, from the Fourier representation (\ref{eq-rep-R}) for
$\mathcal{R}(w(\cdot,t))$, we get the bounds 
\begin{equation}
\label{eq-coeff-decay}
|\alpha_l(\tau)|  \le Ce^{-4\eta\tau}, \qquad  |\beta_l(\tau)|\le Ce^{-4\eta\tau},
\quad \mbox{for all} \,\,  l\ge 1.
\end{equation}
Note that  for $l\ge 3$, we have 
\begin{eqnarray*}
\left |\int_0^s e^{\lambda_l (s-\tau)}\alpha_l(\tau)\cos(l\theta)\, d\tau \right | &\le&
Ce^{\lambda_l s}\int_0^s e^{(-\lambda_l - 4\eta)\tau}\, d\tau \\
&=& \tilde{C}(e^{\lambda_l s} - e^{-4\eta s}) = o(e^{-2s})
\end{eqnarray*}
since $\lambda_l < -2$ and $-4\eta < -2$ for suitably chosen $\eta$. 
This finishes the proof of our claim. 
\end{proof}
By (\ref{eq-semi-rep}), (\ref{eq-zero-part}), (\ref{eq-error-part}) and Claim \ref{claim-error} we get that 
$$w(\theta,\tau) = \left (\alpha_0 + \int_0^s\alpha_0(\tau)e^{-2\tau}\, d\tau \right )\, e^{2s} +  Q_1(\theta,\tau) + Q_2(\theta,\tau)$$
where
\begin{eqnarray*}
Q_1(\theta,\tau) &=& \left (\alpha_1 + \int_0^s\alpha_1(\tau)e^{-\tau}\, d\tau \right  )\, e^s\cos\theta \\
&+&  
\left ( \alpha_2 + \int_0^s\alpha_2(\tau)e^{2\tau}\, d\tau \right )e^{-2s}\cos 2\theta +  o(e^{-2s}),
\end{eqnarray*}
and $Q_2(\theta,\tau)$ is the corresponding expression where the functions 
$\cos \theta$ and $\cos 2\theta $ replaced by  $\sin\theta $ and
$\sin 2\theta$ respectively.   

\begin{claim}We have 
$$\alpha_j + \int_0^s\alpha_j(\tau)e^{-\lambda_j\tau}\, d\tau = 0, \quad \mbox{for} \,\,  j \in \{0,1\}  \quad \mbox{and} \quad 
\beta_1 + \int_0^s\beta_1(\tau)e^{-\lambda_1\tau}\, d\tau = 0.$$
\end{claim}

\begin{proof}
This follows immediately from the representation for
$Q_1(\theta,\tau)$, $Q_2(\theta,\tau)$, the estimate
(\ref{eq-coeff-decay}) and the fact that $|w(\theta,\tau)| \le
Ce^{-2\eta\tau}$,  for some $\eta \in (0,1)$ that can be taken as
close to $1$ as we want. 
\end{proof}

The above discussion implies that 
\begin{equation}
\label{eq-rep-k}
\tilde{\kappa}(\theta,s) = 1 + (a(s)\cos(2\theta) + b(s)\sin(2\theta)) e^{-2s} + o(e^{-2s}),
\end{equation}
where
$$a(s) := \alpha_2 + \int_0^s\alpha_2(\tau)e^{2\tau}\, d\tau \quad  \mbox{and} \quad 
b(s) := \beta_2 + \int_0^s\beta_2(\tau)e^{2\tau}\, d\tau.$$ 

Having the precise asymptotic behavior of $\tilde{\kappa}$ will help us finish the proof of Lemma
\ref{lem-I-zero}. More precisely, we have 
$$\tilde{p}(\theta,\tau) = \tilde{\kappa}^2(\theta,\tau) = 1 + 2\, [a(\tau)\cos 2\theta + b(\tau)\sin 2\theta]\, e^{-2\tau} + o(e^{-2\tau}),$$
which implies
$$\tilde{\alpha}(\theta,\tau) = \tilde{p}_{\theta}(\theta,\tau) = 4\, [b(\tau)\cos 2\theta - a(\tau)\sin 2\theta ] \,  e^{-2\tau}  + o(e^{-2\tau})$$
and 
$$\tilde{\alpha}_{\theta}(\theta,\tau) = 8 \, [-b(\tau)\sin 2\theta - a(\tau)\cos 2\theta ] \, e^{-2\tau} + o(e^{-2\tau}).$$
We can now compute
\begin{eqnarray*}
I(\alpha(t)) &=& e^{4\tau}\int_0^{2\pi}(\tilde{\alpha}_{\theta}^2 - 4\tilde{\alpha}^2)\, d\theta  \\
&=& 64 \, \int_0^{2\pi} [b^2(\tau)\sin^2 2\theta + 2a(\tau)b(\tau)\sin 2\theta\cos 2\theta + a^2(\tau)\cos^22\theta]   \, d\, \theta\\
&-& 64 \, \int_0^{2\pi} [b^2(\tau)\cos^2 2\theta - 2a(\tau)b(\tau)\sin 2\theta\cos 2\theta) +  a^2(\tau)\sin^2 2\theta ] \, d\theta   +   o(1) \\
&=& 64 \, (b(\tau)^2 - a(\tau)^2)\,  \int_0^{2\pi} [\sin^2 2\theta - \cos^2 2\theta]\, d\theta  \\
&+&   256 \, a(\tau)\, b(\tau)\int_0^{2\pi} \sin 2\theta\cos 2\theta\, d\theta + o(1) \\
&=& o(1)
\end{eqnarray*}
since $\int_0^{2\pi}(\sin^2 2\theta - \cos^2 2\theta)\, d\theta = 0$ and
$ \int_0^{2\pi}\sin 2\theta\cos 2\theta\, d\theta = 0$. 
It now follows that
$$\lim_{t\to 0} I(\alpha(t)) = \lim_{\tau\to\infty} e^{4\tau}\int_0^{2\pi} 
(\tilde{\alpha}_{\theta}^2 - 4\tilde{\alpha}^2)\, d\theta = \lim_{\tau\to\infty} o(1) = 0$$
which finishes the proof of Lemma \ref{lem-I-zero}.
\end{proof}

We will now complete the proof  of Theorem \ref{thm-classification}.

\begin{proof}[Proof of Theorem \ref{thm-classification}]
By Proposition \ref{prop-ancient} we have 
$$I(\alpha(t)) \equiv 0, \qquad \mbox{for all} \,\,  t < 0.$$
Lemma \ref{lem-mon-funct} implies that $\alpha_t \equiv 0$, that is,
$$p\, (\alpha_{\theta\theta} + 4\alpha) = 0$$
which means (since $p > 0$) that 
$$\alpha_{\theta\theta} + 4\alpha = 0$$
and therefore 
$$\alpha(\theta,t) = a(t)\cos 2\theta + b(t)\sin 2\theta$$
for some functions in time $a(t)$ and $b(t)$. 
Since $\alpha = p_{\theta}$, by integrating in $\theta$ we conclude that 
\begin{equation}
\label{eq-sol-p}
p(\theta,t) = a(t)\frac{\sin 2\theta}{2} - b(t)\frac{\cos 2\theta}{2} + c(t)
\end{equation}
for another function in time $c(t)$. 
We will now use that $p$ satisfies \eqref{eqn-p} to determine the functions
$a,b$ and $c$. We have:
$$p_t(\theta,t) = a'(t)\frac{\sin 2\theta}{2} -b'(t)\frac{\cos 2\theta}{2} + c'(t)$$
and
$$p_{\theta} = a(t)\cos 2\theta + b(t)\sin 2\theta$$
and
$$p_{\theta\theta} = -2a(t)\sin 2\theta + 2b(t)\cos 2\theta.$$
Plugging those in (\ref{eqn-p}), after a direct  computation,  we obtain
\begin{equation}
\label{eq-abc}
a'(t)\frac{\sin 2\theta}{2} - b'(t)\frac{\cos 2\theta}{2} + c'(t) = -\frac{(a^2(t)+b^2(t))}{2}
+ 2c^2(t)
\end{equation}
which implies that $a'(t) = b'(t) = 0$, that is, $a(t) = a$ and $b(t) = b$ for some constants 
$a,b$,  for all $t < 0$. We can then write our solution $p$ in the form
$$p(\theta,t) = a\frac{\sin 2\theta}{2} + (c(t)-\frac b2)\cos^2\theta + (c(t)+
\frac b2)\sin^2\theta.$$
By the Proposition \ref{prop-limit}, we have 
\begin{equation}
\label{eq-lim-inf1}
\lim_{t\to-\infty}p(\theta,t) = a\frac{\sin 2\theta}{2} + (c-\frac b2)\cos^2\theta + 
(c+\frac b2)\, \sin^2\theta = \lambda\, \cos^2 (\theta+\gamma)
\end{equation}
where $c := \lim_{t\to -\infty} c(t)$, for some fixed angle $\gamma$ and constant $\lambda > 0$.  
Using that $\cos(\theta + \gamma) = \cos\theta\cos\gamma - \sin\theta\sin\gamma$ and (\ref{eq-lim-inf1}) 
we obtain the equation 
\begin{equation*}
\begin{split}
\frac{a}{2}\sin 2\theta  + (c- \frac  b2 -\lambda\cos^2\gamma)\cos^2\theta 
+ (c + \frac b2  &- \lambda\sin^2\gamma)
\sin^2\theta = -\lambda\sin 2\theta\sin\gamma\cos\gamma
\end{split}
\end{equation*}
which is equivalent to
$$\frac{a}{2}\sin 2\theta + c - \frac \lambda2 - (\frac b2 + \frac \lambda 2\cos 2\gamma)\cos 2\theta = -\lambda\sin 2\theta
\sin\gamma\cos\gamma.$$
Hence
\begin{equation}
\label{eq-values}
a = -\lambda\sin 2\gamma, \quad  b = -\lambda\cos 2\gamma, \quad  c = \frac{\lambda}{2}.
\end{equation}
By using the exact expressions for $a$ and $b$ given by
(\ref{eq-values}), the ODE (\ref{eq-abc}) becomes
$$c'(t) = 2\, c^2(t) - \frac{\lambda^2}{2}.$$
We distinguish between two cases.
\begin{case} We have $\lambda = 0$. \\
By (\ref{eq-values}), the expression (\ref{eq-sol-p}) for $p(\theta,t)$ becomes
$$p(\theta,t) = c(t)$$
which satisfies
\begin{equation}
\label{eq-pressure0}
c'(t)= 2\, c^2(t).
\end{equation}
Using that the $\lim_{t\to 0} p(\theta,t) = +\infty$, after integrating (\ref{eq-pressure0})
we get 
$$p(\theta,t) = c(t) = \frac{1}{(-2t)}$$
which corresponds to contracting circles by the curve shortening flow.
\end{case}

\begin{case} We have  $\lambda > 0$. \\
In this case we have
$$p(\theta,t) = \frac{\lambda}{2}\cos 2(\theta+\gamma) + c(t)$$
and we need to solve the ODE
$$\frac{dc}{(2c-\lambda)(2c+\lambda)} = \frac{dt}{2}.$$
Using that the $\lim_{t\to 0} p(\theta,t) = +\infty$ (which implies
the $\lim_{t\to 0} c(t) = +\infty$) yields to 
$$c(t) = \frac{\lambda(1+e^{2\lambda t})}{2(1-e^{2\lambda t})}.$$
Combining the above gives that 
$$p(\theta,t) = \lambda(\frac{1}{1-e^{2\lambda t}} - \sin^2\theta),$$
for a parameter $\lambda > 0$.
\end{case}
\end{proof}
 
\section{Convexity conditions}
\label{sec-equiv}

In this section we will give some  conditions which will 
guarantee that our closed  ancient  embedded solution 
to the CSF is  convex.  For
that we will need the following monotonicity result that has been
already known. We  will  include its proof here  for the
completeness of our exposition. We will also need the Sturm oscillation
theorem (see \cite{CZ}) that says if we have a uniformly elliptic linear equation
$$u_t = a(x,t)u_{xx} + b(x,t)u_x + c(x,t)u,$$
where $a, a_x, a_{xx}, a_t, b, b_t, b_x$ are measurable and bounded
and if $u(x,t)$ is a solution that never vanishes on $\partial I\times
[0,T]$ then the number of zeros of $u(\cdot,t)$ is finite and
non-increasing for all $t\in (0,T)$ and it drops exactly at multiple
zeroes. The set $\{t\in (0,T) : u(\cdot,t)$ has a multiple zero $\}$
is discrete.

\begin{lem}
The total absolute curvature of any curve shortening flow on  closed
curves is decreasing in time.
\end{lem}

\begin{proof}
The proof of the Lemma can be found  in \cite{CZ}. We use the evolution
equation \eqref{eq-strictly-par}. 
For every $\e > 0$, after integration by parts  we have
\begin{eqnarray}
\label{eq-abs-curv}
&&\frac{d}{dt}\int_{\gamma_t}(\e^2+\kappa^2)^{1/2}\, ds \\
&=& \int_{\gamma_t}\kappa\, (\e^2+\kappa^2)^{-1/2} \,  (\kappa_{ss}  + \kappa^3 )\, ds - \int_{\gamma_t}\kappa^2(\e^2+\kappa^2)^{1/2}\, ds  \nonumber \\
&=& -\int_{\gamma_t}\e^2(\e^2+\kappa^2)^{-3/2}\, \kappa_s^2 \, ds -
\int_{\gamma_t} \e^2\kappa^2(\e^2 + \kappa^2)^{-1/2}\, ds \nonumber \\
\qquad \qquad &\le& 0. \nonumber 
\end{eqnarray}
Letting $\e\to 0$ finishes the proof of the Lemma.
\end{proof}  
 
\begin{proof}[Proof of Corollary \ref{cor-equiv}]
Assume that our solution is convex. In that case the Harnack
estimate(\ref{eq-harnack}) implies that $\kappa_t \ge 0$,  which
immediately shows that  $0 < \kappa \le C_1$,  for some uniform constant
$C_1$ and all $t\le t_0 < 0$. From the previous section we already
know that convex solutions have to be either the contracting circles
or the Angenant ovals. In both cases we have
$$\int_{\gamma_t}|\kappa|\, ds \le C_2, \qquad  \mbox{for all} \,\,\, t < 0.$$

Assume now that we have our curvature conditions as in the statement
of our Corollary. We will show the solution needs to be convex.  Take a sequence $t_i\to -\infty$ and $p_i\in S^1$ so that $Q_i
:= |\kappa(p_i,t_i)| = \max|\kappa(\cdot,t_i)|$. Rescale the flow by
setting
$$\gamma_i(p,t) = Q_i(\gamma(p+p_i,tQ_i^{-2}+t_i) - \gamma(p_i,t_i)),
\qquad \mbox{for} \,\, t\in (-\infty,-t_i).$$
 Recall that $Q_i \le C$, uniformly in $i$.
For each $i$, $\gamma_i(\cdot,t)$ solves the curve shortening flow,
$\gamma(0,0) = (0,0)$ and $|\kappa_i| \le C_1$, by our curvature
assumption. We will assume that our solutions $\gamma_i$ are parametrized
by the arc length. We would like to consider a limiting curve that may
not be close so we will assume each $\gamma_i$ is defined on
$\mathbb{R}$ as a periodic map. Since the curvature bound implies that
all its derivatives are uniformly bounded, by the Arzela-Ascoli
theorem we can extract a subsequence which converges to a complete
solution $\tilde{\gamma}$ to (\ref{eq-general}) on every compact
subset of $\mathbb{R}\times (-\infty,\infty)$, so that
$|\tilde{\kappa}(0,0)| = 1$. By the similar arguments to those in
Chapter $5$ in \cite{CZ} one can conclude that $\tilde{\gamma}$ has to
be uniformly convex.  For the sake of completeness we will briefly
sketch the argument from \cite{CZ}.

First notice that any inflection point of $\tilde{\gamma}(\cdot,t)$
must be degenerate.  To see that we argue by contradiction. If not,
there is an inflection point at $(s_0,t_0)$ so that
$|\kappa_t(s_0,t_0)| =: \eta > 0$. By our
convergence, for sufficiently large $i$ and small $\delta > 0$, for each
$t\in [t_0-\delta,t_0]$, there exists $s_i(t)$ near $s_0$ so that
$\kappa_i(s_i(t),t) =0$ and $|(\kappa_i)_s(s_i(t),t)| \ge \eta/2$. One of the terms on the right hand side
of (\ref{eq-abs-curv}) can be estimated as follows
$$-\int_{\gamma_t}\e^2(\e^2+\kappa^2)^{-3/2}\, \kappa_s^2 \, ds 
\le -c_1\int_{-c_2\e}^{c_2\e} \e^2(\e^2 + s^2)^{-\frac{3}{2}}\, ds =
-\frac{2c_1c_2}{(1+c_2^2)^{1/2}}$$ where the positive
constants $c_1,c_2$ depend only on $\eta$. If we integrate
(\ref{eq-abs-curv}) from $t_0-\delta$ to $t_0$ we get 
$$\int_{\gamma_i(t)}(\e^2+\kappa_i^2)^{1/2}\, ds \le \int_{\gamma_i(t_0-\delta)}(\e^2+\kappa^2)^{1/2}\, ds 
- \frac{2\, c_1c_2\, \delta}{(1+c_2^2)^{1/2}}.$$
If we let $\e\to 0$, rewriting the previous estimate for the original
flow,  we get
\begin{equation}
\label{eq-infl-est}
\int_{\gamma(t_i+t_0)}|\kappa|\, ds - \int_{\gamma(t_i+(t_0-\delta))}|\kappa|\, ds \le
\frac{2c_1c_2\delta}{(1+c_2^2)^{1/2}}.
\end{equation}
Since $\int_{\gamma(t)}|\kappa|\, ds$ is decreasing in $t$ and is
uniformly bounded, there is a finite $\lim_{t\to
  -\infty}\int_{\gamma(t)}|\kappa|\, ds$. This implies the left hand
side in (\ref{eq-infl-est}) tends to zero as $i\to \infty$ and we get
a contradiction.

If $\tilde{\gamma}$ were not uniformly convex there would exist some
$t_0$ and $p, q$ so that
$\tilde{\kappa}(p,t_0)\cdot\tilde{\kappa}(q,t_0) < 0$ and therefore
$\tilde{\kappa}(r,t_0) = 0$, for some point $r$ between $p$ and $q$.
By the Sturm oscillation theorem \cite{CZ} applied to
$\tilde{\kappa}(\cdot,t)$, for some $t_1 > t_0$ and $t_1$ close to
$t_0$, $\tilde{\gamma}$ has a non-degenerate inflection point which
contradicts the conclusion in the previous paragraph.

Once we know the blow up limit is uniformly convex we proceed as
follows. If the limiting curve is closed, then for all $i\ge i_0$,
$\gamma_i$ is uniformly convex which implies the same for
$\gamma(\cdot,t_i)$. Since $t_i\to -\infty$ and the curve shortening
flow preserves the convexity, this would mean that $\gamma(\cdot,t)$ is a
convex solution for all $t\in (-\infty,0)$. 

Suppose now that the limiting
curve is complete and noncompact. Assume that there is a $t_1 < 0$ so that
a curve $\gamma(\cdot,t_1)$ has an inflection point and let $l_1$ be
the length of that curve. Take $M > 0$ as big as we want, say so that
$M/C > l_1$,  where $C$ is an  upper bound on $\kappa$. The convergence
to $\tilde{\gamma}$ is uniform on compact subsets, in particular on a
set $(p,t)\in[0,M]\times [-1,1]$. This implies there is some $i_0$ so
that for $i\ge i_0$, curves $\gamma_i(p,t)$ do not have inflection
points for $p\in [0,M]$ (since our limiting curve $\tilde{\gamma}$ is
uniformly convex). In particular, this means a curve
$\gamma(p_{i_0}+p,t_{i_0})$ does not have any inflection points for
$p\in [0,M/C]$, which contains an interval $[0,l_1]$. Considering a
solution $\gamma(p_{i_0}+p,t)$ for $t\ge t_{i_0}$ and $p\in [0,l_1]$,
by the Sturm oscillation theorem the number of zeros of $\kappa(\cdot,t)$
decreases in time and therefore $\gamma(p_{i_0}+p,t_1)$ does not have
any inflection points for $p\in [0,l_1]$. Since the length of
$\gamma(\cdot,t_1)$ is $l_1$, we conclude that the curve
$\gamma(\cdot,t_1)$ can not have any inflection points. In particular,
this implies $\kappa(\cdot,t) > 0$ for all $t < 0$, that is,
$\gamma(\cdot,t)$ is convex for every $t < 0$. The proof of our Corollary
is now complete. 
\end{proof}


\begin{thebibliography}{11}



\bibitem{CZ} Chou,K.-S., Zhu,X.-P., {\em The curve shortening problem}; Chapman \& Hall/CRC, 2001;
ISBN $1-58488-213-1$.

\bibitem{An} Angenent,S. {\em The zero set of a solution of a parabolic equation}; J.Reine Angew.Math.  390 (1988), 79-96.

\bibitem{G}  Gage, M. E. Curve shortening makes convex curves circular. Invent. Math. 76 (1984), no. 2, 357--364. 

\bibitem{GH} Gage,M., Hamilton,R.S. {\em The heat equation shrinking convex plane curves}; J.Diff.Geom. 23 (1996), 69--96.



\bibitem{Gr} Garyson, M. {\em The heat equation shrinks embedded plane
curves to round points}; J. Differential Geom.  26 (1987), no. 2,
285--314.

\bibitem{Gr1} Grayson, M {\em Shortening embedded curves};  The Annals of Mathematics 129 (1989), 71--111.

\bibitem{Ha} Hamilton, R.S. {\em Differential Harnack estimates for parabolic equations}; preprint.

\bibitem{MD} Mitrinovic, D.S. {\em Analytic inequalities}, Springer 1970.

\end{thebibliography}
\end{document}